\documentclass[12pt]{amsart}
\usepackage{amssymb}

\usepackage{amscd}

\newtheorem{theorem}{Theorem}[section]
\newtheorem{corollary}[theorem]{Corollary}

\theoremstyle{definition}

\theoremstyle{remark}

\begin{document}
\title{On rational curves in $n$-space with given normal bundle}
\author{Herbert Clemens}
\date{November, 2000}
\address{Mathematics Department, University of Utah\\
Salt Lake City, UT, 84112, USA}
\email{clemens@math.utah.edu}
\maketitle

\begin{abstract}
The stable rationality of components of the moduli space of (unparametrized)
rational curves in projective $n$-space with fixed normal bundle is proved,
provided these components dominate the moduli space of immersed rational
curves in the plane.
\end{abstract}

\section{Introduction}

In two papers almost twenty years ago, Eisenbud and Van de Ven studied the
variety of parametrized rational curves in $\Bbb{P}^{3}$ and showed that it
is stratified according to the isomorphism classes of the normal bundle of
the curve and that each stratum is rational. They posed the problem of
rational strata in the unparametrized case and rationality was then proved
for ``half'' the cases by Ballico in \cite{B}. Also in 1984 Katsylo \cite{H}
implicitly proved rationality of the full moduli space of unparametrized
rational curves in $\Bbb{P}^{n}$. There seems to have been little attention
to the problem of rationality of strata of rational curves in
higher-dimensional projective spaces, where Eisenbud and Van de Ven said
that the straightforward extension of the their results to higher dimension
was ``to be feared.''

Let $\tilde{S}_{d}^{n}$ denote the space of (unparametrized) immersed
rational curves in $\Bbb{P}^{n}$, that is 
\begin{equation}
\tilde{S}_{d}^{n}=\frac{\left\{ f:\Bbb{P}^{1}\rightarrow \Bbb{P}^{n}:f\
immersive,\ \deg f=d\right\} }{\Bbb{P}GL\left( 2\right) -action\ on\ \Bbb{P}%
^{1}}  \label{0.1}
\end{equation}
The purpose in writing this short note is to draw attention to the
unparametrized case in $\Bbb{P}^{n}$, and in particular to make an
observation that each stratum of the natural stratification of $\tilde{S}%
_{d}^{n}$ by normal-bundle-type has a natural stratification with strata
which are birationally vector bundles over strata of $\tilde{S}_{d}^{2}$.
Using Ballico's result, it is easy to show that $\tilde{S}_{d}^{2}$ is
stably rational, in fact that 
\begin{equation*}
\tilde{S}_{d}^{2}\times \Bbb{A}^{d+1}
\end{equation*}
is rational. This allows us to conclude stably rational for all components $%
S^{\prime }$ of strata of $\tilde{S}_{d}^{n}$ which dominate $\tilde{S}%
_{d}^{2}$, and so to conclude rationality whenever
\begin{equation*}
\dim S^{\prime }-\dim \tilde{S}_{d}^{2}\geq d+1.
\end{equation*}

The author wishes to thank Aaron Bertram, David Eisenbud, Angelo Vistoli and
Igor Dolgechev for their help in understanding this problem.\footnote{%
Partially supported by NSF grant DMS-9970412}

\section{Immersed rational curves}

Let 
\begin{equation*}
f=\left( s^{0},\ldots ,s^{n}\right) :\Bbb{P}^{1}\rightarrow \Bbb{P}^{n}
\end{equation*}
be an immersion. We have an exact sequence 
\begin{equation*}
0\rightarrow \frak{D}_{1}\left( \mathcal{O}_{\Bbb{P}^{1}}\left( 1\right)
\right) \overset{\mu }{\longrightarrow }f^{*}\frak{D}_{1}\left( \mathcal{O}_{%
\Bbb{P}^{n}}\left( 1\right) \right) \rightarrow N_{f}\rightarrow 0
\end{equation*}
where $\frak{D}_{1}\left( L\right) $ is the sheaf of first-order holomorphic
differential operators on sections of the line bundle $L$ and 
\begin{equation*}
\mu \left( a^{i}\frac{\partial }{\partial U^{i}}\right) =a^{i}\frac{\partial
s^{j}}{\partial U^{i}}\frac{\partial }{\partial X^{j}}.
\end{equation*}
Considering $\frak{D}_{1}$ as a left $\mathcal{O}$-module, apply the functor 
\begin{equation*}
Hom\left( \ ,\mathcal{O}_{\Bbb{P}^{1}}\right)
\end{equation*}
to the exact sequence 
\begin{equation*}
0\rightarrow \frak{D}_{1}\left( \mathcal{O}_{\Bbb{P}^{1}}\left( 1\right)
\right) \overset{\mu }{\longrightarrow }f^{*}\frak{D}_{1}\left( \mathcal{O}_{%
\Bbb{P}^{n}}\left( 1\right) \right) \rightarrow N_{f}\rightarrow 0
\end{equation*}
to obtain 
\begin{equation}
0\rightarrow N_{f}^{\vee }\overset{\nu }{\longrightarrow }f^{*}\mathcal{O}_{%
\Bbb{P}^{n}}\left( -1\right) ^{\oplus \left( n+1\right) }\overset{\mu ^{\vee
}}{\longrightarrow }\mathcal{O}_{\Bbb{P}^{1}}\left( -1\right) ^{\oplus
2}\rightarrow 0.  \label{1.1}
\end{equation}

Now a first-order deformation of the map $f$ gives a first-order deformation
of the map $\mu ^{\vee }$ and so an element of 
\begin{equation*}
\mathrm{Hom}\left( N_{f}^{\vee },\mathcal{O}_{\Bbb{P}^{1}}\left( -1\right)
^{\oplus 2}\right) .
\end{equation*}
Thus, by the exact sequence $\left( \ref{1.1}\right) $, we obtain an element 
\begin{equation*}
\upsilon \in \mathrm{Ext}^{1}\left( N_{f}^{\vee },N_{f}^{\vee }\right)
=H^{1}\left( End\left( N_{f}\right) \right)
\end{equation*}
which measures the first-order deformation of the normal bundle $N_{f}$ as a
vector bundle over $\Bbb{P}^{1}$.

By Grothendieck's lemma 
\begin{equation*}
N_{f}=\bigoplus\nolimits_{i=1}^{n-1}\mathcal{O}_{\Bbb{P}^{1}}\left(
a_{i}^{\prime }\right)
\end{equation*}
where, by adjunction, 
\begin{equation*}
\sum\nolimits_{i=1}^{n-1}a_{i}^{\prime }=d\left( n+1\right) -2.
\end{equation*}
Let $L\subseteq \Bbb{P}^{n}$ be a general linear space of codimension $3$.
Then the family of projective spaces 
\begin{equation*}
\overline{Lf\left( x\right) },\ x\in \Bbb{P}^{1},
\end{equation*}
gives a sub-bundle 
\begin{equation*}
\mathcal{O}_{\Bbb{P}^{1}}\left( d\right) ^{n-2}\subseteq N_{f}.
\end{equation*}
Varying the choice of $L,$ these sub-bundles span $N_{f}$. Thus each $%
a_{i}^{\prime }\geq d$ and we put 
\begin{equation*}
a_{i}=a_{i}^{\prime }-d.
\end{equation*}
Then 
\begin{equation*}
N_{f}=\bigoplus\nolimits_{i=1}^{n-1}\mathcal{O}_{\Bbb{P}^{1}}\left(
d+a_{i}\right)
\end{equation*}
where, for all $i$, 
\begin{equation*}
0\leq a_{i}\leq 3d-2
\end{equation*}
and 
\begin{equation*}
\sum\nolimits_{i=1}^{n-1}a_{i}=2d-2.
\end{equation*}

\section{Projection to $\Bbb{P}^{2}$}

Suppose we begin with an immersion 
\begin{equation*}
f_{0}=\left( s^{0},\ldots ,s^{n},0,\ldots ,0,\right) :\Bbb{P}^{1}\rightarrow 
\Bbb{P}^{n+r}.
\end{equation*}
We wish to study the deformations of $f_{0}$ in $\Bbb{P}^{n+r}$. Such a
deformation is given by choosing $t^{k}\in H^{0}\left( \mathcal{O}_{\Bbb{P}
^{1}}\left( d\right) \right) $ for $k=1,\ldots ,r$ and defining 
\begin{equation}
f_{\varepsilon }=\left( f_{0},\varepsilon t^{1},\ldots ,\varepsilon
t^{r}\right) .  \label{1.2}
\end{equation}
We wish to compute the associated element of $H^{1}\left( End\left(
N_{f_{0}}\right) \right) $ in $\left( \ref{1.1}\right) $.

Applying the functor 
\begin{equation*}
R\mathrm{Hom}\left( N_{f_{0}}^{\vee },\ \right)
\end{equation*}
to the sequence $\left( \ref{1.1}\right) $, the infinitesimal deformation $%
N_{f_{\varepsilon }}$ of $N_{f_{0}}$ is given by the image of the matrix 
\begin{equation}
\left( 
\begin{array}{llllll}
0 & \ldots & 0 & \frac{\partial t^{1}}{\partial U^{0}} & \ldots & \frac{
\partial t^{r}}{\partial U^{0}} \\ 
0 & \ldots & 0 & \frac{\partial t^{1}}{\partial U^{1}} & \ldots & \frac{
\partial t^{r}}{\partial U^{1}}
\end{array}
\right) \in \mathrm{Hom}\left( N_{f_{0}}^{\vee },\mathcal{O}_{\Bbb{P}%
^{1}}\left( -1\right) ^{\oplus 2}\right)  \label{1.9}
\end{equation}
in $\mathrm{Ext}^{1}\left( N_{f_{0}}^{\vee },N_{f_{0}}^{\vee }\right) =%
\mathrm{Ext}^{1}\left( N_{f_{0}},N_{f_{0}}\right) $, that is, by the
equivalence class of $\left( \ref{1.9}\right) $ in 
\begin{equation*}
\frac{\mathrm{Hom}\left( N_{f_{0}}^{\vee },\mathcal{O}_{\Bbb{P}^{1}}\left(
-1\right) ^{\oplus 2}\right) }{\mathrm{Hom}\left( N_{f_{0}}^{\vee },f^{*}%
\mathcal{O}_{\Bbb{P}^{n}}\left( -1\right) ^{\oplus \left( n+r+1\right)
}\right) \circ \mu ^{\vee }}.
\end{equation*}

Now 
\begin{equation*}
N_{f_{0}}^{\vee }=\mathcal{O}_{\Bbb{P}^{1}}\left( -d\right) ^{r}\oplus
N_{g}^{\vee }
\end{equation*}
where 
\begin{equation*}
g=\left( s^{0},\ldots ,s^{n}\right) :\Bbb{P}^{1}\rightarrow \Bbb{P}^{n}
\end{equation*}
and the projection 
\begin{equation*}
\mathrm{Hom}\left( N_{f_{0}}^{\vee },\mathcal{O}_{\Bbb{P}^{1}}\left(
-1\right) ^{\oplus 2}\right) \rightarrow \mathrm{Hom}\left( \mathcal{O}_{%
\Bbb{P}^{1}}\left( -d\right) ^{r},\mathcal{O}_{\Bbb{P}^{1}}\left( -1\right)
^{\oplus 2}\right)
\end{equation*}
takes $\mathrm{Hom}\left( N_{f_{0}}^{\vee },f^{*}\mathcal{O}_{\Bbb{P}%
^{n}}\left( -1\right) ^{\oplus \left( n+r+1\right) }\right) $ to the
homomorphisms generated by the columns of 
\begin{equation*}
\left( 
\begin{array}{lll}
\frac{\partial s^{0}}{\partial U^{0}} & \ldots & \frac{\partial s^{n}}{
\partial U^{0}} \\ 
\frac{\partial s^{0}}{\partial U^{1}} & \ldots & \frac{\partial s^{n}}{
\partial U^{1}}
\end{array}
\right) .
\end{equation*}
Thus, for example, every first-order deformation of $N_{f_{0}}$ which does
not change the normal bundle is a combination of a deformation of $g$ in $%
\Bbb{P}^{n}$ combined with a deformation of $\Bbb{P}^{n}$ in $\Bbb{P}^{n+r}$.

Finally now suppose that $n=2$. Then we have that 
\begin{equation*}
N_{f_{0}}^{\vee }=\mathcal{O}_{\Bbb{P}^{1}}\left( -d\right) ^{r}\oplus
N_{g}^{\vee }=\mathcal{O}_{\Bbb{P}^{1}}\left( -d\right) ^{\oplus r}\oplus 
\mathcal{O}_{\Bbb{P}^{1}}\left( 2-3d\right) .
\end{equation*}
Also 
\begin{equation*}
Hom\left( N_{f_{0}}^{\vee },\ \right) =Hom\left( \mathcal{O}_{\Bbb{P}
^{1}}\left( -d\right) ^{r},\ \right) \oplus Hom\left( N_{g}^{\vee },\
\right) 
\end{equation*}
and 
\begin{eqnarray*}
\mathrm{Ext}^{1}\left( N_{f_{0}}^{\vee },N_{f_{0}}^{\vee }\right)  &=&%
\mathrm{Ext}^{1}\left( \mathcal{O}_{\Bbb{P}^{1}}\left( -d\right) ^{\oplus
r},N_{f_{0}}^{\vee }\right)  \\
&=&\mathrm{Ext}^{1}\left( \mathcal{O}_{\Bbb{P}^{1}}\left( -d\right)
,N_{g}^{\vee }\right) ^{\oplus r} \\
&=&H^{1}\left( \mathcal{O}_{\Bbb{P}^{1}}\left( 2-2d\right) \right) ^{\oplus
r}.
\end{eqnarray*}

\section{The case of $\Bbb{P}^{2+1}\label{P3case}$}

From what we have just seen in the last section, the $j$-th entry in the
element of $\mathrm{Ext}^{1}\left( N_{f_{0}}^{\vee },N_{f_{0}}^{\vee
}\right) $ corresponding to the deformation $\left( \ref{1.2}\right) $ is
given by applying 
\begin{equation*}
R\mathrm{Hom}\left( \mathcal{O}_{\Bbb{P}^{1}}\left( -d\right) ,\ \right)
\end{equation*}
to the sequence 
\begin{equation*}
0\rightarrow N_{g}^{\vee }\rightarrow \mathcal{O}_{\Bbb{P}^{1}}\left(
-d\right) ^{\oplus 3}\rightarrow \mathcal{O}_{\Bbb{P}^{1}}\left( -1\right)
^{\oplus 2}\rightarrow 0
\end{equation*}
and finding the image $\upsilon $ of 
\begin{equation}
\left( 
\begin{array}{l}
\frac{\partial t}{\partial U^{0}} \\ 
\frac{\partial t}{\partial U^{1}}
\end{array}
\right) \in \mathrm{Hom}\left( \mathcal{O}_{\Bbb{P}^{1}}\left( -d\right) ,%
\mathcal{O}_{\Bbb{P}^{1}}\left( -1\right) ^{\oplus 2}\right)  \label{new2}
\end{equation}
in 
\begin{equation*}
\frac{\mathrm{Hom}\left( \mathcal{O}_{\Bbb{P}^{1}}\left( -d\right) ,\mathcal{%
O}_{\Bbb{P}^{1}}\left( -1\right) ^{\oplus 2}\right) }{\left\{ \Bbb{C} 
\begin{array}{l}
\frac{\partial s^{0}}{\partial U^{0}} \\ 
\frac{\partial s^{0}}{\partial U^{1}}
\end{array}
+\Bbb{C} 
\begin{array}{l}
\frac{\partial s^{1}}{\partial U^{0}} \\ 
\frac{\partial s^{1}}{\partial U^{1}}
\end{array}
+\Bbb{C} 
\begin{array}{l}
\frac{\partial s^{2}}{\partial U^{0}} \\ 
\frac{\partial s^{2}}{\partial U^{1}}
\end{array}
\right\} }=H^{1}\left( \mathcal{O}_{\Bbb{P}^{1}}\left( 2-2d\right) \right) .
\end{equation*}

On the other hand, for all integers $a$ we have 
\begin{eqnarray*}
\mathrm{Ext}^{1}\left( N_{f_{0}}^{\vee },N_{f_{0}}^{\vee }\right) &=&\mathrm{%
\ Ext}^{1}\left( N_{f_{0}}^{\vee }\left( a\right) ,N_{f}^{\vee }\left(
a\right) \right) \\
&=&H^{1}\left( End\left( N_{f_{0}}^{\vee }\left( a\right) \right) \right)
\end{eqnarray*}
and so we have a pairing 
\begin{equation}
H^{0}\left( N_{f_{0}}^{\vee }\left( a\right) \right) \otimes H^{1}\left(
End\left( N_{f_{0}}^{\vee }\left( a\right) \right) \right) \rightarrow
H^{1}\left( N_{f_{0}}^{\vee }\left( a\right) \right)  \label{new1}
\end{equation}
which measures the obstruction to deforming sections of $N_{f_{0}}^{\vee
}\left( a\right) $ to first order with a first-order deformation of $%
N_{f_{0}}^{\vee }$ given by an element of 
\begin{equation*}
H^{1}\left( End\left( N_{f_{0}}^{\vee }\right) \right) =H^{1}\left(
End\left( N_{f_{0}}^{\vee }\left( a\right) \right) \right) .
\end{equation*}
But we are in a situation in which the only obstructions to deforming
sections of $N_{f_{0}}^{\vee }\left( a\right) $ to sections of $%
N_{f_{\varepsilon }}^{\vee }\left( a\right) $ for all $\varepsilon $ are of
first order. To see this use $\left( \ref{1.1}\right) $ to write 
\begin{equation*}
0\rightarrow N_{f_{\varepsilon }}^{\vee }\left( a\right) \overset{\nu }{
\longrightarrow }\mathcal{O}_{\Bbb{P}^{n}}\left( a-d\right) ^{\oplus 4}%
\overset{\mu ^{\vee }}{\longrightarrow }\mathcal{O}_{\Bbb{P}^{1}}\left(
a-1\right) ^{\oplus 2}\rightarrow 0
\end{equation*}
with the matrix $\mu ^{\vee }$ is given by 
\begin{equation*}
\left( 
\begin{array}{cccc}
\frac{\partial s^{0}}{\partial U^{0}} & \frac{\partial s^{1}}{\partial U^{0}}
& \frac{\partial s^{2}}{\partial U^{0}} & \varepsilon \frac{\partial t}{%
\partial U^{0}} \\ 
\frac{\partial s^{0}}{\partial U^{1}} & \frac{\partial s^{1}}{\partial U^{1}}
& \frac{\partial s^{2}}{\partial U^{1}} & \varepsilon \frac{\partial t}{%
\partial U^{1}}
\end{array}
\right) .
\end{equation*}
So the equations for extension of sections become 
\begin{equation*}
\left( 
\begin{array}{cccc}
\frac{\partial s^{0}}{\partial U^{0}} & \frac{\partial s^{1}}{\partial U^{0}}
& \frac{\partial s^{2}}{\partial U^{0}} & \varepsilon \frac{\partial t}{%
\partial U^{0}} \\ 
\frac{\partial s^{0}}{\partial U^{1}} & \frac{\partial s^{1}}{\partial U^{1}}
& \frac{\partial s^{2}}{\partial U^{1}} & \varepsilon \frac{\partial t}{%
\partial U^{1}}
\end{array}
\right) \left( 
\begin{array}{c}
\alpha _{00}+\varepsilon \alpha _{01} \\ 
\alpha _{10}+\varepsilon \alpha _{11} \\ 
\alpha _{20}+\varepsilon \alpha _{21} \\ 
\alpha _{30}
\end{array}
\right) =\left( 
\begin{array}{c}
0 \\ 
0
\end{array}
\right)
\end{equation*}
which have no terms of degree $>1$ in $\varepsilon .$

Now we rewrite the pairing $\left( \ref{new1}\right) $ vertically as 
\begin{equation*}
\begin{array}{c}
\left( H^{0}\left( \mathcal{O}_{\Bbb{P}^{1}}\left( a-d\right) \right) \oplus
H^{0}\left( \mathcal{O}_{\Bbb{P}^{1}}\left( a+2-3d\right) \right) \right)
\otimes H^{1}\left( \mathcal{O}_{\Bbb{P}^{1}}\left( 2-2d\right) \right) \\ 
\downarrow \\ 
H^{1}\left( \mathcal{O}_{\Bbb{P}^{1}}\left( a-d\right) \right) \oplus
H^{1}\left( \mathcal{O}_{\Bbb{P}^{1}}\left( a+2-3d\right) \right)
\end{array}
.
\end{equation*}
Thus for $d-1\leq a\leq 3d-2$ the pairing $\left( \ref{new1}\right) $
becomes the multiplication map 
\begin{equation}
\begin{array}{c}
H^{0}\left( \mathcal{O}_{\Bbb{P}^{1}}\left( a-d\right) \right) \otimes
H^{1}\left( \mathcal{O}_{\Bbb{P}^{1}}\left( 2-2d\right) \right) \\ 
\downarrow \\ 
H^{1}\left( \mathcal{O}_{\Bbb{P}^{1}}\left( a+2-3d\right) \right)
\end{array}
.  \label{1.6}
\end{equation}
which we restrict to the element 
\begin{equation*}
image\left( 
\begin{array}{l}
\frac{\partial t}{\partial U^{0}} \\ 
\frac{\partial t}{\partial U^{1}}
\end{array}
\right) \in H^{1}\left( \mathcal{O}_{\Bbb{P}^{1}}\left( 2-2d\right) \right) .
\end{equation*}

It will be important to note one structural property of the map $\left( \ref
{1.6}\right) $. Let $p\in \Bbb{P}^{1}$ with local coordinate $z.$ Using the
exact sequence 
\begin{equation*}
0\rightarrow \mathcal{O}_{\Bbb{P}^{1}}\left( \left( 2-2d\right) \cdot
p\right) \rightarrow \mathcal{O}_{\Bbb{P}^{1}}\rightarrow \left\{ \frac{%
\sum\nolimits_{j=1}^{\infty }\Bbb{C\cdot }z^{j}}{z^{2d-2}}\right\}
\rightarrow 0
\end{equation*}
we rewrite $\left( \ref{1.6}\right) $ as 
\begin{equation}
\begin{array}{c}
\left\{ \sum\nolimits_{i=0}^{a-d}\Bbb{C\cdot }z^{i}\right\} \otimes \left\{ 
\frac{\sum\nolimits_{j=1}^{\infty }\Bbb{C\cdot }z^{j}}{z^{2d-2}}\right\} \\ 
\downarrow \\ 
\left\{ \frac{\sum\nolimits_{j=1}^{\infty }\Bbb{C\cdot }z^{j}}{z^{3d-a-2}}
\right\}
\end{array}
\label{1.7}
\end{equation}
and make the important remark that the dimension of the kernel of any map 
\begin{equation*}
\xi :\left\{ \sum\nolimits_{i=0}^{a-d}\Bbb{C\cdot }z^{i}\right\} \rightarrow
\left\{ \frac{\sum\nolimits_{j=1}^{\infty }\Bbb{C\cdot }z^{j}}{z^{3d-a-2}}
\right\}
\end{equation*}
given by multiplication by $\xi \in \left\{ \frac{\sum\nolimits_{j=1}^{%
\infty }\Bbb{C\cdot }z^{j}}{z^{2d-2}}\right\} $ depends only on the degree
of the leading term in $\xi $. Also, since the map 
\begin{eqnarray*}
H^{0}\left( \mathcal{O}_{\Bbb{P}^{1}}\left( d\right) \right) &\rightarrow
&H^{1}\left( \mathcal{O}_{\Bbb{P}^{1}}\left( 2-2d\right) \right) \\
t &\mapsto &\xi
\end{eqnarray*}
is linear, the leading coefficient must vanish on a codimension-one subspace
of $H^{0}\left( \mathcal{O}_{\Bbb{P}^{1}}\left( d\right) \right) $.

\section{The theorem for $\Bbb{P}^{2+r}$}

So, returning to the general situation of $n=2$ and $r$ arbitrary, we can
completely characterize the extension 
\begin{equation*}
N_{f_{\varepsilon }}^{\vee }
\end{equation*}
by knowing 
\begin{equation*}
h^{0}\left( N_{f_{\varepsilon }}^{\vee }\left( a\right) \right)
\end{equation*}
for each integer $a$. But, for a given first-order deformation 
\begin{equation*}
\Xi \in H^{1}\left( End\left( N_{f_{0}}\right) \right) =H^{1}\left(
End\left( N_{f_{0}}\left( a\right) \right) \right) ,
\end{equation*}
these latter numbers are given by the dimensions of the subspaces of $%
H^{0}\left( N_{f_{0}}^{\vee }\left( a\right) \right) $ consisting of
sections which deform to first-order, and hence to all orders, with $%
\varepsilon $. But these are just the dimension of 
\begin{equation*}
\ker \left( H^{0}\left( N_{f_{0}}\left( a\right) \right) \rightarrow
H^{1}\left( N_{f_{0}}\left( a\right) \right) \right) .
\end{equation*}
So these dimensions are given by the kernels of the $r$ maps $\left( \ref
{1.7}\right) .$ Now use the map 
\begin{equation*}
\begin{array}{c}
\Phi :H^{0}\left( \mathcal{O}_{\Bbb{P}^{1}}\left( d\right) \right) ^{\oplus
r}\rightarrow \left\{ \frac{\sum\nolimits_{j=1}^{\infty }\Bbb{C\cdot }z^{j}}{
z^{2d-2}}\right\} ^{\oplus r} \\ 
\left( t^{1},\ldots ,t^{r}\right) \mapsto \left( 
\begin{array}{lll}
\frac{\partial t^{1}}{\partial U^{0}} & \ldots & \frac{\partial t^{r}}{%
\partial U^{0}} \\ 
\frac{\partial t^{j}}{\partial U^{1}} & \ldots & \frac{\partial t^{j}}{%
\partial U^{1}}
\end{array}
\right)
\end{array}
\end{equation*}
to define a map 
\begin{equation*}
H^{0}\left( \mathcal{O}_{\Bbb{P}^{1}}\left( d\right) \right) ^{\oplus
r}\rightarrow \left\{ 1,\ldots ,2d-3\right\} ^{\oplus r}
\end{equation*}
which associates to each $\left( t^{1},\ldots ,t^{r}\right) $ the degrees of
the leading coefficients of $\Phi \left( t^{1},\ldots ,t^{r}\right) $. This
is a semi-continuous map for which the preimage of each element is a linear
space minus a linear subspace. And as we vary the immersed curve $g\left( 
\Bbb{P}^{1}\right) \subseteq \Bbb{P}^{2}$, the dimensions of each linear
space is locally constant on a Zariski open set. Thus:

\begin{theorem}
\label{only}i) Let $S$ be the set of all smooth rational curves in $\Bbb{P}%
^{2+r}$ whose image via the standard projection 
\begin{equation*}
\Bbb{P}^{2+r}\dashrightarrow \Bbb{P}^{2}
\end{equation*}
is a fixed immersed curve and whose normal bundle is 
\begin{equation}
\bigoplus\nolimits_{i=1}^{r+1}\mathcal{O}_{\Bbb{P}^{1}}\left( d+a_{i}\right) 
\label{100}
\end{equation}
for some fixed value of $\left( a_{1},\ldots ,a_{r+1}\right) .$ Then $S$ is
rational.

ii) Suppose the general curve in $\tilde{S}_{d}^{2}$ deforms to a curve in $%
\tilde{S}_{d;a_{1},\ldots ,a_{r+1}}^{r+2}$. The set 
\begin{equation*}
\tilde{S}_{d;a_{1},\ldots ,a_{r+1}}^{r+2}
\end{equation*}
consisting of all smooth rational curves in $\Bbb{P}^{2+r}$ with normal
bundle $\left( \ref{100}\right) $ is has a component which is birationally
isomorphic to a vector bundle over the set of immersed rational curves in $%
\Bbb{P}^{2}$, and all other components consists of curves projecting into
proper subvarieties of $\tilde{S}_{d}^{2}$.
\end{theorem}

\begin{proof}
We consider the set 
\begin{equation*}
\tilde{g}\cdot GL\left( 2\right)
\end{equation*}
of reparametrizations of a fixed map 
\begin{equation*}
\tilde{g}:\Bbb{A}^{2}\rightarrow \Bbb{A}^{3}
\end{equation*}
with projectivization 
\begin{equation*}
g:\Bbb{P}^{1}\rightarrow \Bbb{P}^{2}.
\end{equation*}
We need only check the action of the group $GL\left( 2\right) $ on the space 
\begin{equation}
\left( \tilde{g}\cdot GL\left( 2\right) \right) \times H^{0}\left( \mathcal{O%
}_{\Bbb{P}^{1}}\left( d\right) \right) ^{\oplus r}  \label{5.1}
\end{equation}
parametrization of maps 
\begin{equation*}
\tilde{f}:\Bbb{A}^{2}\rightarrow \Bbb{A}^{3+r}.
\end{equation*}
But, considering $\left( \ref{5.1}\right) $ as a (trivial) vector bundle
over the affine variety $\left( \tilde{g}\cdot GL\left( 2\right) \right) $,
this is a free action of the group 
\begin{equation*}
\frac{GL\left( 2\right) }{\mu _{d}I}
\end{equation*}
on the vector bundle $\left( \ref{5.1}\right) $ over a free action on the
base space $\left( \tilde{g}\cdot GL\left( 2\right) \right) $ and this group
acts as a group of vector bundle isomorphisms. So, by descent theory for
coherent sheaves and faithful flatness, the quotient is a vector bundle with
fiber isomorphic to $H^{0}\left( \mathcal{O}_{\Bbb{P}^{1}}\left( d\right)
\right) ^{\oplus r}$.
\end{proof}

\begin{corollary}
Suppose the general curve in $\tilde{S}_{d}^{2}$ deforms to a curve in $%
\tilde{S}_{d;a_{1},\ldots ,a_{r+1}}^{r+2}$. Then each component $S^{\prime }$
of $\tilde{S}_{d;a_{1},\ldots ,a_{r+1}}^{r+2}$ over the generic curve of $%
\tilde{S}_{d}^{2}$ is stably rational. Furthermore $S^{\prime }$ is rational
if
\begin{equation*}
\left( \dim S^{\prime }-\dim \tilde{S}_{d}^{2}\right) \geq \left( d+1\right)
.
\end{equation*}
\end{corollary}

\begin{proof}
Katsylo proved \cite{H} that 
\begin{equation*}
\tilde{S}_{d}^{3}=\overline{\tilde{S}_{d;d-1,d-1}^{3}}
\end{equation*}
is rational. But by Theorem \ref{only} $\tilde{S}_{d;d-1,d-1}^{3}$ is
birationally a vector bundle over $\tilde{S}_{d}^{2}$. Since vector bundles
are locally trivial in the\thinspace Zariski topology, we conclude that the
product of $\tilde{S}_{d}^{2}$ with the vector-bundle fiber $F\cong \Bbb{A}%
^{d+1}$ is rational. That is
\begin{equation*}
\tilde{S}_{d;d-1,d-1}^{3}\overset{birat.}{\longleftrightarrow }\tilde{S}%
_{d}^{2}\times \Bbb{A}^{d+1}.
\end{equation*}
So $\tilde{S}_{d}^{2}$ is stably rational. Now suppose a component $%
S^{\prime }$ of $\tilde{S}_{d;a_{1},\ldots ,a_{r+1}}^{r+2}$ dominates $%
\tilde{S}_{d}^{2}$ via the map
\begin{equation*}
\tilde{S}_{d;a_{1},\ldots ,a_{r+1}}^{r+2}\dashrightarrow \tilde{S}_{d}^{2}
\end{equation*}
induced from the projection map
\begin{equation*}
\Bbb{P}^{r+2}\dashrightarrow \Bbb{P}^{2}.
\end{equation*}
Again by Theorem \ref{only}, $S^{\prime }$ is birationally a (locally
trivial) vector bundle over $\tilde{S}_{d}^{2}$ with fiber $F^{\prime }\cong 
\Bbb{A}^{s}$. So $S^{\prime }$ is always stably rational and is rational
whenever
\begin{equation*}
s\geq d+1.
\end{equation*}
\end{proof}

Notice that some of the $\tilde{S}_{d;a_{1},\ldots ,a_{r+1}}^{r+2}$ may not
dominate $\tilde{S}_{d}^{2}$, in which case we can conclude nothing.

\end{document}